\documentclass[12pt,leqno,draft]{article}
\usepackage{amsfonts}
\usepackage{amsmath, dsfont, amsthm, amsfonts, amssymb, color}
\usepackage{mathrsfs}
\usepackage{color}
\usepackage{stmaryrd}
\newtheorem{thm}{Theorem}[section]

\newtheorem{lem}[thm]{Lemma}

\theoremstyle{definition}

\newcommand{\scr}[1]{\mathscr #1}
\definecolor{wco}{rgb}{0.5,0.2,0.3}

\numberwithin{equation}{section} \theoremstyle{remark}

\newcommand{\ua}{\uparrow}

\title{{\bf
  Singular McKean-Vlasov (Reflecting) SDEs with Distribution Dependent Noise}\footnote{Supported in
 part by  NNSFC (11831014, 11801406, 11921001). } }
\author{
{\bf   Xing Huang $^{a)}$,  Feng-Yu Wang $^{a), b)}$  }\\
\footnotesize{ a)Center for Applied Mathematics, Tianjin
University, Tianjin 300072, China}\\
\footnotesize{  xinghuang@tju.edu.cn}\\
 \footnotesize{ b)Department of Mathematics, Swansea University, Bay Campus, SA1 8EN, United Kingdom}\\
\footnotesize{  wangfy@tju.edu.cn}}
\begin{document}
\allowdisplaybreaks
\def\R{\mathbb R}  \def\ff{\frac} \def\ss{\sqrt} \def\B{\mathbf
B} \def\W{\mathbb W}
\def\N{\mathbb N} \def\kk{\kappa} \def\m{{\bf m}}
\def\ee{\varepsilon}\def\ddd{D^*}
\def\dd{\delta} \def\DD{\Delta} \def\vv{\varepsilon} \def\rr{\rho}
\def\<{\langle} \def\>{\rangle} \def\GG{\Gamma} \def\gg{\gamma}
  \def\nn{\nabla} \def\pp{\partial} \def\E{\mathbb E}
\def\d{\text{\rm{d}}} \def\bb{\beta} \def\aa{\alpha} \def\D{\scr D}
  \def\si{\sigma} \def\ess{\text{\rm{ess}}}
\def\beg{\begin} \def\beq{\begin{equation}}  \def\F{\scr F}
\def\Ric{\text{\rm{Ric}}} \def\Hess{\text{\rm{Hess}}}
\def\e{\text{\rm{e}}} \def\ua{\underline a} \def\OO{\Omega}  \def\oo{\omega}
 \def\tt{\tilde} \def\Ric{\text{\rm{Ric}}}
\def\cut{\text{\rm{cut}}} \def\P{\mathbb P} \def\ifn{I_n(f^{\bigotimes n})}
\def\C{\scr C}      \def\aaa{\mathbf{r}}     \def\r{r}
\def\gap{\text{\rm{gap}}} \def\prr{\pi_{{\bf m},\varrho}}  \def\r{\mathbf r}
\def\Z{\mathbb Z} \def\vrr{\varrho} \def\ll{\lambda}
\def\L{\scr L}\def\Tt{\tt} \def\TT{\tt}\def\II{\mathbb I}
\def\i{{\rm in}}\def\Sect{{\rm Sect}}  \def\H{\mathbb H}
\def\M{\scr M}\def\Q{\mathbb Q} \def\texto{\text{o}} \def\LL{\Lambda}
\def\Rank{{\rm Rank}} \def\B{\scr B} \def\i{{\rm i}} \def\HR{\hat{\R}^d}
\def\to{\rightarrow}\def\l{\ell}\def\iint{\int}
\def\EE{\scr E}\def\Cut{{\rm Cut}}
\def\A{\scr A} \def\Lip{{\rm Lip}}
\def\BB{\scr B}\def\Ent{{\rm Ent}}\def\L{\scr L}
\def\R{\mathbb R}  \def\ff{\frac} \def\ss{\sqrt} \def\B{\mathbf
B}
\def\N{\mathbb N} \def\kk{\kappa} \def\m{{\bf m}}
\def\dd{\delta} \def\DD{\Delta} \def\vv{\varepsilon} \def\rr{\rho}
\def\<{\langle} \def\>{\rangle} \def\GG{\Gamma} \def\gg{\gamma}
  \def\nn{\nabla} \def\pp{\partial} \def\E{\mathbb E}
\def\d{\text{\rm{d}}} \def\bb{\beta} \def\aa{\alpha} \def\D{\scr D}
  \def\si{\sigma} \def\ess{\text{\rm{ess}}}
\def\beg{\begin} \def\beq{\begin{equation}}  \def\F{\scr F}
\def\Ric{\text{\rm{Ric}}} \def\Hess{\text{\rm{Hess}}}
\def\e{\text{\rm{e}}} \def\ua{\underline a} \def\OO{\Omega}  \def\oo{\omega}
 \def\tt{\tilde} \def\Ric{\text{\rm{Ric}}}
\def\cut{\text{\rm{cut}}} \def\P{\mathbb P} \def\ifn{I_n(f^{\bigotimes n})}
\def\C{\scr C}      \def\aaa{\mathbf{r}}     \def\r{r}
\def\gap{\text{\rm{gap}}} \def\prr{\pi_{{\bf m},\varrho}}  \def\r{\mathbf r}
\def\Z{\mathbb Z} \def\vrr{\varrho} \def\ll{\lambda}
\def\L{\scr L}\def\Tt{\tt} \def\TT{\tt}\def\II{\mathbb I}
\def\i{{\rm in}}\def\Sect{{\rm Sect}}  \def\H{\mathbb H}
\def\M{\scr M}\def\Q{\mathbb Q} \def\texto{\text{o}} \def\LL{\Lambda}
\def\Rank{{\rm Rank}} \def\B{\scr B} \def\i{{\rm i}} \def\HR{\hat{\R}^d}
\def\to{\rightarrow}\def\l{\ell}
\def\8{\infty}\def\I{1}\def\U{\scr U} \def\n{{\mathbf n}}
\maketitle

\begin{abstract}   By using Zvonkin's transformation and a two-step fixed point argument in distributions,  the well-posedness and regularity estimates are derived for singular  McKean-Vlasov  SDEs with   distribution dependent noise, where  the drift  contains a term growing linearly in space and distribution and a locally integrable term independent of distribution, while  the noise coefficient is weakly differentiable in space and Lipschitz continuous in   distribution with respect to the sum of Wasserstein and weighted variation distances.   The main results  extend  existing ones derived  for noise coefficients either independent of distribution, or having nice linear functional derivatives in distribution.  Singular reflecting  SDEs with distribution dependent noise are also studied. \end{abstract}

 \vskip 2cm

\section{Introduction}

As a crucial stochastic model characterizing nonlinear Fokker-Planck equations and mean field particle systems,   the following McKean-Vlasov (i.e. distribution dependent) SDE has  been intensively investigated:
  \beq\label{E1} \d X_t= b_t(X_t , \L_{X_t})\d t+\si_t(X_t,\L_{X_t})\d W_{t},\ \
t\in [0,T],\end{equation}
where $T>0$ is a fixed constant,
  $(W_t)_{t\in [0,T]}$ is an $m$-dimensional Brownian motion on a complete filtration probability space $(\OO,\{\F_t\}_{t\in [0,T]},\P)$, $\L_{X_t}$ is the law of $X_t$, and for the space $\scr P$ of
  probability measures on $\R^d$ equipped with the weak topology,
\begin{align*} &b: [0,T]\times\R^d\times \scr P \to \R^d,\ \ \si: [0,T]\times \R^d\times \scr P \to \R^d\otimes\R^m
 \end{align*}
 are measurable. Among many other references,     see for instance \cite{BB,BBP,CR,CF,HW,HW20,L,MV,RZ,FYW1,ZG}.

When the noise coefficient $\si_t(x,\mu)=\si_t(x)$,  by using Zvonkin's  transformation,  the well-posedness, regularity estimates and exponential ergodicity have been studied in \cite{21Ren, W21c, W21d}  for the drift $b_t(x,\mu)$ containing  a time-spatial locally integrable term in $\tt L_{p}^q(T)$ for some $(p,q)\in \scr K$ introduced in \cite{XXZZ}, see \eqref{KK0} and \eqref{KK1} below.

Concerning singular McKean-Vlasov SDEs, the well-posedness is derived in \cite{CF, ZG} when
the noise coefficient $\si_t(x,\mu)$ has a nice linear functional derivative in $\mu$ besides other conditions, where in \cite{CF} the drift $b_t(x,\mu)$  is bounded and uniformly Lipschitz continuous in $\mu$ with respect to the total variation distance, and in \cite{ZG} the drift $b_t(x,\mu)$ is Lipschitz continuous in $\mu$ with respect to a weighted variation distance uniformly in $(t,x)$, and
$ \|b_\cdot(\cdot,\mu)\|_{\tt L_p^q(T)}<\infty$ uniformly in $\mu$ for some $(p,q)\in \scr K$.

Comparing with \cite{CF, ZG},    this paper studies \eqref{E1} for $\si_t(x,\cdot)$ not necessarily having   linear functional derivatives,   and for    $b_t(x,\mu)$  unbounded in $\mu$ and containing a singular distribution independent term.  For instance, let $\si_t(x,\mu)=\si(\mu):= f(\mu)   I_{d\times d}$, where $k\ge 1,$
$I_{d\times d}$ is the identity matrix, and $f(\mu):= 1+ \mu(|\cdot|^k) \land 1. $ Then $\si$ is  Lipschitz continuous in
the $k^{th}$-Wasserstein distance  and hence satisfies  assumption $(A_1)$ introduced below, but it does not have bounded linear  functional derivative required in \cite{CF, ZG} , according to (2.3) in \cite{CF} and the fact that $f$ is not Lipschitz continuous in the total variation norm.

Instead of the usual fixed point method developed for the well-posedness of distribution dependent SDEs, we will adopt a two-step fixed point argument by freezing the distribution variables in $b$ and $\si$ respectively.

Let  $k\in [1,\infty)$. Then
$$\scr P_k=\bigg\{\mu\in \scr P: \|\mu\|_k:=\mu(|\cdot|^k)^{\ff 1 k}:=\left(\int_{\R^d}|x|^k\mu(\d x)\right)^{\frac{1}{k}}<\infty\bigg\} $$ is a Polish space under the $k^{th}$-Wasserstein distance $\W_k$:
  $$\W_k(\mu, \nu):= \inf_{\pi\in \C(\mu,\nu)} \bigg(\int_{\R^d\times\R^d} |x-y|^k \pi(\d x,\d y)\bigg)^{\ff 1 {k}},\ \ \mu,\nu\in \scr P_k,$$
  where $\C(\mu,\nu)$ is the set of all couplings of $\mu$ and $\nu$. Moreover, $\scr P_k$ is a complete metric space under the weighted variation norm
  $$\|\mu-\nu\|_{k, var} := \sup_{|f|\le 1+|\cdot|^k} \big|\mu(f)-\nu(f)\big|,\ \ \ \mu,\nu\in \scr P_k.$$
  By \cite[Theorem 6.15]{Villani}, there exists a constant $\kk>0$ such that
  \beq\label{GPP} \|\mu-\nu\|_{var}+\W_k(\mu,\nu)^k\le \kk \|\mu-\nu\|_{k,var},\end{equation}
  where $\|\cdot\|_{var}$ is the total variation norm. On the other hand, when $k>1$
  there is no any constant $c>0$ such that $\|\mu-\nu\|_{k,var}\ge c \W_k(\mu,\nu)$ holds for all $\mu,\nu\in \scr P_k$.

We call equation \eqref{E1} strongly (weakly) well-posed for distributions in $\scr P_k$, if for any $\F_0$-measurable initial value
$X_0$ with $\L_{X_0}\in \scr P_k$ (respectively any initial distribution $\mu\in \scr P_k$), it has a unique strong solution (respectively weak solution) such that $\L_{X_\cdot}\in C([0,T];\scr P_k)$, the space of continuous maps from $[0,T]$ to   the Polish space  $(\scr P_k,\W_k)$. Moreover, we call \eqref{E1} well-posed for distributions in $\scr P_k$ if it is strongly  and weakly well-posed for distributions in $\scr P_k$. 
In this case, we denote
$$P_t^\ast\mu=\L_{X_t}\ \ \text{for\ the\ solution\ with} \ \L_{X_0}=\mu\in \scr P_k.$$

 To measure the   singularity  of $b_t(x,\mu)$ in $(t,x)$, we recall locally integrable functional spaces introduced in \cite{XXZZ}.
For any $t>s\ge 0$ and $p,q\in (1,\infty)$, we write $f\in  \tilde{L}_p^q([s,t])$  if $f: [s,t]\times\R^d\to \R$ is measurable with
 $$\|f\|_{\tilde{L}_p^q([s,t])}:= \sup_{z\in\R^d}\bigg\{\int_{s}^t \bigg(\int_{B(z,1)}|f(u,x)|^p\d x\bigg)^{\ff q p} \d u\bigg\}^{\ff 1q}<\infty,$$ where $B(z,1):=\{x\in \R^d: |x-z|\le 1\}$ is the unit ball centered at point $z$.
 When $s=0$, we simply denote
\beq\label{KK0} \tt L_p^q(t)=\tt L_p^q([0,t]),\ \ \|f\|_{\tilde{L}_p^q(t)}=\|f\|_{\tilde{L}_p^q([0,t])}.\end{equation}
 We will take $(p,q)$ from the   space
 \beq\label{KK1} \scr K:=\Big\{(p,q): p,q>2, \ff{d}p+\ff 2 q<1\Big\}.\end{equation}
 For any $\mu\in C([0,T];\scr P_k)$,  let
 $$\si_t^\mu(x):= \si_t(x,\mu_t),\ \ \ b_t^\mu(x):= b_t(x,\mu_t),\ \ \ (t,x)\in [0,T]\times \R^d.$$
 We make the following assumption.

\beg{enumerate} \item[$(A_0)$] There exist   constants $K>K_0\ge 0$, $l\in\mathbb N$, $\{(p_i,q_i): 0\le i\le l\}\subset \scr K$ and
$1\le f_i\in \tt L^{q_i}_{p_i}(T)$ for $0\le i\le l$ such that
$\si^\mu_t(x)$ and $b^\mu_t(x):= b_t^{(1)}(x)+ b_t^{\mu,0}(x)$ satisfy the following conditions for all $\mu\in C([0,T];\scr P_k)$.
\item[$(1)$] $a^\mu:= \si^\mu(\si^\mu)^*$ is invertible  with
$\|a^\mu\|_\infty +\|(a^\mu)^{-1}\|_\infty\le K$ and
$$\lim_{\vv\downarrow 0}\sup_{\mu\in C([0,T];\scr P_k)} \sup_{t\in [0,T], |x-y|\le \vv} \|a^\mu_t(x)-a^\mu_t(y)\|= 0.$$
\item[$(2)$] $b^{(1)}$ is locally bounded on $[0,T]\times \R^d$, $\si^\mu_t$ is weakly differentiable such that
\beg{align*}& |b_t^{\mu,0}(x)|\le f_0(t,x)+K_0\|\mu_t\|_k,\ \
  \|\nn \si_t^\mu(x)\|\le \sum_{i=1}^l f_i(t,x),\ \ \ (t,x)\in [0,T]\times\R^d,\\
  &|b_t^{(1)}(x)- b_t^{(1)}(y)|\le K |x-y|,\ \ t\in [0,T], x,y\in \R^d.\end{align*} \end{enumerate}
This assumption implies the well-posedness of the SDE with drift $b_t^\mu(x)$ and noise coefficient $\si_t^\nu(x)$  for all $\mu,\nu\in C([0,T]; \scr P_k)$, see
\cite[Theorem 2.1]{21Ren}.  To prove the well-posedness of   \eqref{E1}, we need the following conditions
on the distribution dependence.

\beg{enumerate} \item[$(A_1)$]  For any $t\in [0,T], x\in \R^d$ and $\mu,\nu\in \scr P_k$,
$$\|\si_t(x,\mu)-\si_t(x,\nu)\|+ |b_t(x,\mu)-b_t(x,\nu)|\le \W_k(\mu,\nu) \sum_{i=0}^l f_i(t,x).$$
  \end{enumerate}

Our first result is the following.

\beg{thm}\label{T1} Assume $(A_0)$ and $(A_1)$. Then the following assertions hold.
 \begin{enumerate}
 \item[$(1)$]   $\eqref{E1}$ is well-posed for   distributions in $\scr P_{k}$. Moreover, for any $j\ge k$ there exists a constant $c(j)>0$ such that    the solution satisfies   \beq\label{NES} \E\Big[ \sup_{t\in [0,T]} |X_t|^{j} \big|\F_0\Big] \le c(j)\big\{1+|X_0|^j +  (\E[|X_0|^k])^{\ff j k}\big\}.\end{equation}
\item[$(2)$] For any $N>0$ and $j\ge k$, there exists a constant $C_{j,N}>0$ such that for any two solutions $X_t^i$ of $\eqref{E1}$ with
$\E[|X_0^i|^k]\le N, i=1,2$,
\beq\label{NWS} \E\Big(\sup_{t\in [0,T]} |X_t^1-X_t^2|^j\Big|\F_0\Big)\le C_{j,N} \big\{|X_0^1-X_0^2|^j+ (\E[|X_0^1-X_0^2|^k] )^{\ff j k}\big\}.\end{equation}Consequently,
  \beq\label{NWS'} \sup_{t\in [0,T]} \W_k(P_t^*\mu^1, P_t^*\mu^2)\le 2 C_{k,N} \W_k(\mu^1, \mu^2),\ \ \mu^1,\mu^2\in \scr P_k, \ \mu^1(|\cdot|^k),\mu^2(|\cdot|^k)\le N.\end{equation}
When $K_0=0$, this estimate holds for some constant $C_j>0$ replacing $C_{j, N}$ for any two solutions for distributions in $\scr P_k$.
\end{enumerate}\end{thm}

Comparing with $(A_1)$, the following assumption  allows  weaker distribution dependence for $b_t(x,\cdot)$ but needs  $b^{(1)}=0$ and stronger conditions on $\si$.

\

\beg{enumerate} \item[$(A_2)$]   $b^{(1)}=0$, and there exists a constant $\kk\ge 0$ such that the following conditions hold for all $t\in [0,T], x,y\in \R^d$ and $\mu,\nu\in \scr P_k$.
 \beg{align*} &|b_t(x,\mu)-b_t(x,\nu)|\le  \big\{\kk \|\mu-\nu\|_{k,var}+\W_{k}(\mu,\nu)\big\}\sum_{i=0}^l f_i(t,x),\\
  &\|\si_t(x,\mu)\|^2\lor \|(\si_t\si_t^*)^{-1}(x,\mu)\|\le K,\\
&\|\si_t(x,\mu)-\si_t(y,\nu)\|\le K\big(|x-y|+\W_k(\mu, \nu)\big),\\
&\|\{\si_t(x,\mu)-\si_t(y,\mu)\}-\{\si_t(x,\nu)-\si_t(y,\nu)\}\|\le K|x-y|\W_k(\mu, \nu).
\end{align*} \end{enumerate}

\paragraph{Remark 1.1.}   It is easy to see that the fourth inequality in $(A_2)$ holds if $\si_t(x,\mu)$ is differentiable in $x$ with
$$\|\nn \si_t(\cdot,\mu)(x)- \nn \si_t(\cdot,\nu)(x)\|\le K \W_k(\mu,\nu),\ \ \mu,\nu\in \scr P_k, x\in\R^d.$$
Indeed, this implies
\beg{align*} &\|\{\si_t(x,\mu)-\si_t(y,\mu)\}-\{\si_t(x,\nu)-\si_t(y,\nu)\}\|\\
&= \bigg\|\int_0^1 \big\{\nn_{x-y} \si_t(y+s(x-y), \mu) - \nn_{x-y} \si_t(y+s(x-y), \nu)\big\}\d s\bigg\|\\
&\le \int_0^1 \big\|\nn_{x-y} \si_t(y+s(x-y), \mu) - \nn_{x-y} \si_t(y+s(x-y), \nu)\big\| \d s\\
&\le K|x-y|\W_k(\mu, \nu).\end{align*}

  \beg{thm}\label{T2} Assume $(A_0)$ and $(A_2)$. Then Theorem $\ref{T1}(1)$ holds.
   If $\kk=0$, then for any $N\geq1$, there exists a constant $C(N)>0$, such that
\beq\label{NES2} \|P_t^\ast\mu-P_t^\ast\nu\|_{var}\leq \frac{C(N)}{\sqrt{t}}\W_{k}(\mu,\nu),\ \ t>0,\|\mu\|_k\vee\|\nu\|_k\leq N.\end{equation}
If moreover $K_0=0$, then the constant $C(N)$ can be independent of $N$.

\end{thm}

The above two theorems are proved in Sections 2 and 3 respectively,  and Theorem \ref{T1} will be extended in Section 4 to reflecting SDEs.

\section{Proof of Theorem \ref{T1}}

Let us  explain the main idea of the two-step fixed point argument.

Let $X_0$ be $\F_0$-measurable with $\gg:= \L_{X_0}\in \scr P_k$. Let
$$\C_k^\gg:=\big\{\mu\in C([0,T];\scr P_k):\ \mu_0=\gg\big\}.$$
We   solve \eqref{E1} with a fixed  distribution parameter   $\mu\in \C_k^\gg$ in the drift:
\beq\label{ENN}   \d X_t^{\mu}= b_t(X_t^{\mu}, \mu_t)\d t+\si_t(X_t^{\mu},\L_{X_t^{\mu}})\d W_{t},\ \
t\in [0,T], X_0^{\mu}=X_0,\end{equation}
such that the well-posedness of this SDE for distributions in $\scr P_k$ provides a map
$$\C_k^\gg\ni\mu\mapsto   \Phi^\gg_\cdot\mu:= \L_{X_\cdot^{\mu}} \in\C_k^\gg.$$
Then   the well-posedness of \eqref{E1} follows if  the map $\Phi^\gg$ has a unique fixed point in $\C_k^\gg$.

To solve   \eqref{ENN}, we further fix the   distribution parameter $\nu\in \C_k^\gg$ in $\si$
such that the SDE becomes
$$  \d X_t^{\mu,\nu} = b_t(X_t^{\mu,\nu}, \mu_t)\d t+\si_t(X_t^{\mu,\nu},\nu_t)\d W_{t},\ \
t\in [0,T],  X_0^{\mu,\nu}=X_0,$$     which is well-posed under $(A_0)$ according to \cite[Theorem 2.1]{21Ren}. This gives a map
\beq\label{EMN2}  \C_k^\gg \ni\nu\mapsto \Phi_\cdot^{\gg,\mu}\nu:= \L_{X_\cdot^{\mu,\nu}}\in \C_k^\gg.\end{equation}
So, we  first prove that this map has a unique fixed point such that    \eqref{ENN} is well-posed, then   apply the fixed point theorem to  $\Phi^\gg$ to derive the well-posedness of the original   SDE \eqref{E1}.

For any $\kk\ge 0$, let
$$\W_{k,\kk var}(\mu^1,\mu^2):= \W_k(\mu^1,\mu^2)+ \kk \|\mu^1-\mu^2\|_{k,var},\ \ \mu^1,\mu^2\in \scr P_k.$$
  To apply the fixed point theorem, we will use the following complete metrics on $\C_k^\gg$ for $\theta>0$ and $\kk\ge 0$:
 \beg{align}\label{DIS} \nonumber& \W_{k,\kk var, \theta}(\mu,\nu):= \sup_{t\in [0,T]} \e^{-\theta t} \W_{k,\kk var} (\mu_t,\nu_t),\\
 & \W_{k, \theta}(\mu,\nu):= \sup_{t\in [0,T]} \e^{-\theta t} \W_{k} (\mu_t,\nu_t),\ \ \mu,\nu\in \scr C_k^\gg.\end{align}

To  prove that $\Phi^\gg$  has a unique fixed point in $\C_k^\gg$, we need to restrict the map to the following     bounded
    subspaces of $\C_k^\gg$:
   \beq\label{KGN} \C_k^{\gg,N} := \Big\{\mu\in \C_k^\gg: \sup_{t\in [0,T]} \e^{-Nt} (1+\mu_t(|\cdot|^k)) \le N\Big\},\ \ N>0, \end{equation}
   and to prove that these spaces are $\Phi^\gg$-invariant for large $N$. This enables us to verify the contraction of $\Phi^\gg$ in $\C_k^{\gg,N}$ under a suitable complete metric.

For this purpose, we present the following lemmas.
 The first one ensures the well-posedness of \eqref{ENN}.

    \beg{lem} \label{LN1} Assume  $(A_0)$ and that for some constant $\kk\ge 0$,
    \beq\label{NN*1}\beg{split} & |b_t(x,\nu_1)-b_t(x,\nu_2)|\le \W_{k,\kk var}(\nu_1,\nu_2)\sum_{i=0}^l f_i(t,x),\\
&\|  \si_t(x,\nu_1)-\si_t(x,\nu_2)\|\le  \W_k(\nu_1,\nu_2)\sum_{i=0}^l f_i(t,x)\end{split}\end{equation} holds for   any $\nu_1,\nu_2 \in \scr P_k, t\in [0,T]$ and $x\in\R^d$.
Then
  $\eqref{ENN}$ is well-posed for distributions in    $\scr P_k$. Moreover,  there exist $\theta_0>0$ and decreasing function  $\bb: [\theta_0,\infty)\to (0,\infty)$ with $\beta(\theta)\downarrow 0$ as $\theta\uparrow \infty$ such that
\beq \label{WDI}  \W_{k,\theta}(\Phi^\gamma \mu,\Phi^\gamma \nu)\leq \bb(\theta)  \W_{k,\kk var,\theta}(\mu,\nu),\ \   \mu,\nu\in \C_k^{\gg,N}.
\end{equation}

\end{lem}

\beg{proof}
(a) For the well-posedness,  it suffices to prove that   $\Phi^{\gg,\mu}$  defined in \eqref{EMN2} has a unique fixed point in  $\C_k^\gg$.

In general, let      $\mu^i \in\C_k^{\gg^i,N}$  for some    $N>0, \gg^i\in \scr P^k, i=1,2.$
For   $\nu^i \in \C_k^{\gg^i}$ and initial value $X_0^i$ with $\L_{X_0^i}=\gg^i, i=1,2,$ consider the SDEs
\beq\label{W*1}  \d X_{t}^{i}= b_t^{\mu^i}(X_t^{i})\d t+\si_t^{\nu^i}(X_t^{i})\d W_{t},\  \  t\in [0,T], i=1,2. \end{equation}
According to \cite[Theorem 2.1]{21Ren}, under   $(A_0)$  these SDEs are well-posed, and by \cite[Theorem 2.1]{YZ},  there exist  constants  $c_0,\ll_0\ge  0$  depending on   $N$  via   $\mu^1\in\C_k^{\gg,N}$  due to
  $$|b_t^{\mu^1,0}(x)|\le f_0(t,x) + K_0 \|\mu_t^1\|_k,$$
such that for any   $\ll\ge \ll_0$,  the  PDE
  \beq\label{W*2}
\Big(\pp_t +\ff 12 {\rm tr}\{a_t^{\nu^1}\nn^2\} \Big) u_t+(\nn u_t)b_t^{\mu^1} =  \ll u_t- b_t^{\mu^1,0}, \ \ t\in [0,T], u_T=0\end{equation}
has a unique solution
such that
   \beq\label{3AA} \|\nn^2u\|_{\tt L_{p_0}^{q_0}(T)} \le c_0,\ \ \|u\|_\infty+\|\nn u\|_\infty\le \ff 1 2.\end{equation}
Let $Y_t^i:= \Theta_t(X_t^{i}), i=1,2, \Theta_t:=id +u_t$.
By It\^o's formula we obtain
{  \beg{align*}
  &\d Y_t^1=\big\{b_t^{(1)} + \lambda u_t \big\} (X_{t}^{1}) \d t+ (\{\nabla\Theta_t\}\sigma_t^{\nu^1} )(X_{t}^{1})\,\d W_t,\\
 & \d Y_t^2 =\big\{\big\{b_t^{(1)} + \lambda u_t+{  (\nabla\Theta_t)(b^{\mu^2}_t-b^{\mu^1}_t)\big\}(X_{t}^{2})} \\
 &  + {  \ff 1 2\big[{\rm tr} \{(a_t^{\nu^2}- a_t^{\nu^1}) \nn^2   u_t\}\big](X_{t}^{2})} \big\}\d t+(\{\nabla\Theta_t\}\sigma_t^{\nu^2})(X_{t}^{2})\,\d W_t.\end{align*}}
Let $\eta_{t}:= |X_{t}^{1}- X_{t}^{2}|$ and
   \beg{align*} &g_r:=  \sum_{i=0}^l f_i  (r, X_{r}^{2}),\ \ \tt g_r:= g_r  \|\nn^2  u_r(X_{r}^{2})\|,\\
   &\bar g_r:= \sum_{i=1}^2 \|\nn^2 u_r\|(X_r^i)+\sum_{j=1}^2 \sum_{i=0}^l f_i  (r, X_{r}^{j}),\ \ r\in [0,T].\end{align*}
 Since $b_t^{(1)} + \lambda u_t $ is Lipschitz continuous uniformly in $t\in [0,T]$, by  $(A_0)$, \eqref{NN*1} and the maximal functional inequality in \cite[Lemma 2.1]{XXZZ},
     there exists a constant
 $c_1>0$ depending on $N$ such that
   \beg{align*}&\big| \big\{b_r^{(1)} + \lambda u_r \big\} (X_{r}^{1})  - \big\{b_r^{(1)} + \lambda u_r \big\} (X_{r}^{2}) \big|\le c_1 \eta_r,\\
 & \big| \big\{(\nabla\Theta_r)(b^{\mu^2}_r-b^{\mu^1}_r)\big\}(X_{r}^{2})\big|\le c_1  g_r  \W_{k,\kk var} (\mu^1_r,\mu^2_r),\\
 &\big|\big[{\rm tr} \{(a_r^{\nu^2}- a_r^{\nu^1}) \nn^2   u_r\}\big](X_{r}^{2})\big|\le c_1 \tt g_r \W_k(\nu_r^1,\nu_r^2),\\
 &\big\| \big\{(\nabla\Theta_r ) \sigma_r^{\nu^1} \big\} (X_{r}^{1})- \big\{(\nabla\Theta_r ) \sigma_r^{\nu^2} \big\} (X_{r}^{2})\big\|\\
 &\le c_1  \bar g_r \eta_r+  c_1 g_r \W_{k}(\nu^1_r,\nu^2_r),\ \ r\in [0,T].\end{align*}
So, by It\^o's formula, for any $j\ge k$
   we find a constant    $c_2>1$  depending on  $N$  such that
\beq\label{**1}  \d |Y_t^1-Y_t^2|^{2j}\le c_2  \eta_t^{2j} \d A_t + c_2 (g_t^2+\tt g_t) \big\{\W_{k, \kk var}(\mu^1_t,\mu_t^2)^{2j} + \W_k(\nu^1_t,\nu_t^2)^{2j} \big\}\d t +\d M_t\end{equation}
   holds for some martingale $M_t$   with $M_0=0$ and
   $$A_t:=\int_0^t \big\{1+g_s^2+\tt g_s+ \bar g_s^2\big\}\d s.$$
   Since $\|\nn u\|_\infty\le \ff 1 2$ implies $|Y_t^1-Y_t^2|\ge \ff 1 2 \eta_t,$ this implies
     \beq\label{NNP} \beg{split}& \eta_{t}^{2j}  \le 2^{2j}M_t+ 2^{2j} \eta_0^{2j} +  2^{2j}c_2\int_0^t \eta_{r}^{2j} \d A_r \\
     &+ 2^{2j}c_2 \int_0^t  (g_s^2+\tt g_s) \big\{\W_{k,\kk var}(\mu^1_s,\mu_s^2)^{2j} + \W_k(\nu^1_s,\nu_s^2)^{2j} \big\}\d s  \end{split}\end{equation}
  for some constant $c_2>0$ and all $t\in [0,T]$.
 By   \eqref{3AA},  $f_i\in \tt L^{q_i}_{p_i}(T)$ for $(p_i,q_i)\in \scr K$,   Krylov's and  Khasminskii's estimates (see \cite{YZ}), we find an increasing function
 $\aa: (0,\infty)\to (0,\infty)$ and a decreasing function $\vv: (0,\infty)\to (0,\infty)$ with $\vv_\theta\to 0$ as $\theta\to\infty$, such that
  {  $$\E[\e^{r A_T}|\F_0]\le \aa(r),\ \ \ r>0,$$}
{  $$  \sup_{t\in [0,T]} \E\bigg( \int_0^t \e^{-2k \theta(t-r)}(g_r^2+\tt g_r)\d r \bigg|\F_0\bigg)\le \vv_\theta, \ \ \theta>0.$$}
  By  the stochastic Gronwall inequality
and  the maximal inequality (see \cite{XXZZ}), we find   a   constant  {  $c_3>0$}  depending on {  $N$}  such that {  \eqref{NNP} } yields
\beq\label{*F}\beg{split}   &\Big\{\E\Big(\sup_{s\in [0,t]}  \eta_{s}^j\Big|\F_0\Big)\Big\}^{2} \\
& \le c_3 \E \bigg(\eta_0^{2j}+ \int_0^t  (g_s^2+\tt g_s) \big\{\W_{k,\kk var}(\mu^1_s,\mu_s^2)^{2j} + \W_k(\nu^1_s,\nu_s^2)^{2j} \big\}\d s\bigg|\F_0\bigg) \\
 &\le c_3  \eta_0^{2j}+c_3  \e^{2k\theta t}   \vv_\theta     \big\{\W_{k,\kk var,\theta}(\mu^1,\mu^2)^{2j} +\W_{k,\theta}(\nu^1,\nu^2)^{2j}  \big\}.\end{split}\end{equation}
  Noting that
  {  $$  \W_k(\L_{X_{t}^{1}},\L_{X_{t}^{2}})^k \le  \E[|X_t^1-X_t^2|^k]=  \E [\eta_{t}^k],$$} by taking $j=k$ we obtain
{   \beq\label{W} \W_{k,\theta} (\L_{X^{1}},\L_{X^2})^k
      \le   \ss{c_3} \E[\eta_0^k]+  \ss{c_3\vv_\theta } \big\{\W_{k,\kk var,\theta}(\mu^1,\mu^2)^k+ \W_{k,\theta}(\nu^1,\nu^2)^k\big\}.\end{equation}}
By taking $X_0^1=X_0^2=X_0$ and    $\mu^1=\mu^2=\mu\in \C_k^{\gg,N}$, 
when    $\theta>0$ is large enough such that   $\sqrt{c_3\vv_\theta}\le \ff 1 2$,   $\Phi^{\gg,\mu}\nu^i=\L_{X^i}$  satisfies
{  $$\W_{k,\theta}(\Phi^{\gg,\mu}\nu^1, \Phi^{\gg,\mu}\nu^2)\le \ff 1 2 \W_{k,\theta}(\nu^1,\nu^2),\ \ \nu_1,\nu_2\in \C_k^\gg.$$}
Thus,  $\Phi^{\gg,\mu}$  has a unique fixed point in $\C_k^\gg$, so that \eqref{ENN} is well-posed for distributions in $\scr P_k$.

(b) Taking $\nu^i= \Phi^\gg \mu^i$, we have $\L_{X^i}= \Phi^\gg\mu^i$, so that \eqref{W} becomes
$$ \W_{k,\theta} (\Phi^\gg\mu^1,\Phi^\gg\mu^2)
      \le    (c_3\vv_\theta)^{\ff 1 {2k}} \big\{\W_{k,\kk var,\theta}(\mu^1,\mu^2)+ \W_{k,\theta} (\Phi^\gg\mu^1,\Phi^\gg\mu^2) \big\}.$$ Taking $\theta_0>0$ large enough such that
      $c_3\vv_{\theta_0} <1$  we prove \eqref{WDI} for
      $$\bb(\theta):= \ff { (c_3\vv_\theta)^{\ff 1 {2k}} }{1- (c_3\vv_\theta)^{\ff 1 {2k}} }, \ \ \theta\ge \theta_0.$$

\end{proof}



\beg{lem}\label{LN2}  Assume {  $(A_0)$}.\beg{enumerate}
   \item[$(1)$] There exists a constant {  $N_0>0$} such that for any {  $N\geq N_0$}    we have
   {      $\Phi^\gg\C_k^{\gg,N}\subset  \C_k^{\gg,N}$}.
   \item[$(2)$] Solutions to {  $\eqref{E1}$}  for distributions in {  $\scr P_k$ } satisfy  {  $\eqref{NES}$} for any $j\ge k$ and some constant $c(j)>0$.
   \end{enumerate}
   \end{lem}

   \beg{proof}
{  (1)} Simply denote {   $  M_t = \int_0^t \si_s(X_s^\mu,\L_{X_s^\mu})\d W_s.$}
Since $\|\si\|_\infty<\infty$ due to  {  $(A_0)$},  we have
$$\sup_{t\in [0,T]} \E[|M_t|^k]<\infty.$$
Combining this with
Lemma \ref{LJJ}  below,  we find some constants {  $c_0,c_1 >0$} such that
   \beg{align*}
 &\E(1+|X_t^\mu|^k) \\
 &\le   \E(1+|X_0|^k) + c_0\E\left|\int_0^t ( K_0 \|\mu_s\|_k+ f_0(s,X_s^\mu)+
|X_s^\mu|+1) \d s\right|^k
 + \E\big|M_t\big|^k\\
 &\le c_1+c_1\left|\int_0^t  \|\mu_s\|^2_k \d s\right|^{k/2}+c_1 \int_0^t \E(1+|X_s^\mu|^k)\d s, \ \ t\in [0,T].\end{align*}
 By Gronwall's inequality, we find $c_2,c_3>0$ such that
 \beg{align*}    &\E(1+|X_t^\mu|^k)   \leq c_2+c_2\left|\int_0^t {  \e^{- \frac{2N}{k}s}\|\mu_s\|^2_k} \e^{\frac{2N}{k}s}\d s\right|^{k/2}\\
 &\leq c_3+c_3N^{1-k/2}\e^{Nt},\ \ \mu\in\C_k^{\gg,N},  t\in [0,T]. \end{align*}
Therefore, we   find   a constant {  $N_0>0$}  such that
{  $$\sup_{t\in [0,T]} (1+\|\Phi_{t}^\gg\mu\|^k_k) \e^{-Nt} \le c_3+c_3N^{1-k/2}\leq N,\ \ N\ge N_0, \mu\in\C_k^{\gg,N}.$$}
That is, {  $\Phi^\gg\C_k^{\gg,N}\subset\C_k^{\gg,N}$} for {  $N\ge N_0$}.

{  (2) }   Let {  $X_t$}  solve  \eqref{E1}   with   $\gg:=\L_{X_0}\in\scr P_k$, and denote
{  $\mu_t:=\L_{X_t }$. } Then {  $X_t=X_t^{\mu}$}. By {  $(A_0)$}  and It\^o's formula,
for any $j\ge 1$ we find a constant {  $c_1>0$ } such that
\beq\label{**3}     |X_t|^{2j} -|X_0|^{2j}
 \le    c_1 \int_0^t \big\{ 1+ |X_s|^{2j} + |X_s|^{2j-1}f_0(s,X_s)+ \|\mu_s\|_k^{2j}  \big\}\d s + M_t\end{equation}
holds for some martingale {  $M_t$}  with {  $\d\<M\>_t\le c_1^2 |X_t|^{2(2j-1)}\d t.$ }
Noting that
{  \beg{align*} &c_1\int_0^t   |X_s|^{2j-1}f_0(s,X_s)\d s\le c_1\Big(\sup_{s\in [0,t]}  |X_s|^{2j-1} \Big) \int_0^t   f_0(s,X_s)\d s \\
&\le \ff 1 2 \sup_{s\in [0,t]} |X_s|^{2j}  + c_2\bigg(\int_0^tf_0(s,X_s)\d s\bigg)^{2j} \end{align*} }
holds for some constant  $c_2>0$, we see that {  $\eta_t:= \sup_{s\in [0,t]} |X_s|^{2j} $}
 satisfies
   {  \beq\label{NN*3} \beg{split} \eta_t\le 2|X_0|^{2j}+&\,2 c_1 \int_0^t \big\{1+\eta_s+\|\mu_s\|_k^{2j}\big\}\d s +2 c_2\bigg(\int_0^t  f_0(s,X_s)\d s\bigg)^{2j}
  + 2\sup_{s\in [0,t]}M_s.\end{split} \end{equation} }
By {  $\d\<M\>_t\le c_1^2 |X_t|^{2(2j-1)} \d t$ }  and   BDG's inequality,
we find   constants $c_3,c_4>0$  such that
{  \beg{align*} & \E\Big(\sup_{s\in [0,t]} M_s\Big|\F_0\Big) \le c_3 \E\bigg[\bigg(\int_0^t |X_s|^{2(2j-1)}\d s\bigg)^{\ff 1 2}\bigg|\F_0\bigg]\\
&\le \ff 1 4 \E\big(\eta_t\big|\F_0\big)+   c_4\int_0^t\big\{1+\E(\eta_s|\F_0)\big\}\d s.\end{align*} }
Combining this with {  \eqref{NN*3}}  and  { \eqref{JJ2}} below,   we find a constant {  $c_5>0$}  such that
\beq\label{**4}  \E\big(\eta_t\big|\F_0\big)\le c_5 + c_5 |X_0|^{2j}+c_5 \int_0^t \big\{\E(\eta_s|\F_0)+ \|\mu_s\|_k^{2j}\big\}\d s,\ \ t\in [0,T].\end{equation}
 By Gronwall's inequality, there exists a constant $c_6>0$ such that
\beq\label{NN*20} \E\big(\eta_t\big|\F_0\big)\le c_6 + c_6 |X_0|^{2j}+c_6 \int_0^t \|\mu_s\|_k^{2j}\d s,\ \ t\in [0,T].\end{equation}
In particular, choosing {  $j=k$ } and applying Jensen's inequality,   we derive
\beg{align*} &\E\Big[\sup_{s\in [0,t]}|X_s|^{k}\Big|\F_0\Big]\le \Big\{\E\big(\eta_t\big|\F_0\big)\Big\}^{\ff 1 2} \\
&\le \ss{c_6}\big(1+|X_0|^k\big)+\ff{c_6}2 \int_0^t \|\mu_s\|_k^k\d s +\ff 1 2 \sup_{s\in [0,t]} \|\mu_s\|_k^k.\end{align*}
Noting that $\|\mu_s\|_k^k=\E[|X_s|^k]$, by taking expectation we obtain
$$\|\mu_t\|_k^k \le\E\Big[ \sup_{s\in [0,t]} |X_s|^k\Big]\le 2 \ss{c_6}\big(1+\E[|X_0|^k]\big)+c_6 \int_0^t  \|\mu_s\|_k^k\d s,\ \ t\in [0,T].$$
By  Gronwall's inequality, we find a constant $c>0$ such that
{  $$\|\mu_t\|_k^k \le c (1+\E[|X_0|^k]),\ \ t\in [0,T].$$}
Substituting into {  \eqref{NN*20}} we prove {  $\eqref{NES}$}. \end{proof}

\begin{lem}\label{LJJ}   Assume $(A_0)$.  For any $(p,q)\in \scr K$, there exist a constant $c_0\ge 1 $ and a  function $c:[1,\infty)\to(0,\infty)$ such that for any $j\geq 1$    and $\mu\in \C_k^\gg$, the solution to $\eqref{ENN}$ satisfies
\beq\label{JJ1}
 \E\big[\e^{\int_0^t | f_s(X_{s}^\mu)|^2\d s}\big|\F_0\big] \le  \e^{c_0+c_0 \int_0^t\|\mu_s\|_k^2\d s+c_0 \|f\|_{\tt L^q_p(t)}^{c_0}},   \end{equation}
\beq\label{JJ2}  \E\bigg[\left(\int_0^t |f_s(X_{s}^{\mu})|^2\d s\right)^j\bigg|\F_0\bigg]\leq c(j) \left(1+ \int_0^t\|\mu_s\|_k^2\d s \right)^j \|f\|_{\tt L^{q}_p(t)}^{2j}\end{equation}
for any $t\in [0,T]$ and $   f\in \tt L^q_p(t), t\in [0,T]. $

\end{lem}

\begin{proof}
Consider the SDE
$$\d \bar X_{t}=b_t^{(1)}(\bar X_t)\d t + \sigma_t(\bar X_{t}, \Phi^\gamma_{t}\mu)\d W_t,\ \ \bar X_0=X_0, t\in[0,T].$$
By Khasminskii's estimate (see \cite{YZ}), there exists a constant $c_1>1$   such that
\beq\label{NN9}
 \E\big[\e^{\int_0^t | f_s(\bar X_{s}^\mu)|^2\d s}\big|\F_0\big] \le  \e^{c_1+  c_1 \|f\|_{\tt L^q_p(t)}^{c_1}},  \ \   f\in \tt L_q^p(t), t\in [0,T].\end{equation}
By  $(A_0)$,
$$\xi_t:= \sigma_t(\bar X_{t}, \Phi^\gamma_{t}\mu)^*\{\sigma_t(\bar X_{t}, \Phi^\gamma_{t}\mu)\sigma_t(\bar X_{t}, \Phi^\gamma_{t}\mu)^*\}^{-1} b_t^{\mu,0}(\bar X_t)$$
satisfies
$$|\xi_t|\le c_2 f_0(t, \bar X_t) +c_2\|\mu_t\|_k,\ \ t\in [0,T]$$
for some constant $c_2>0$. Combining this with \eqref{NN9},  we conclude that
$$R_t:= \e^{\int_0^t\<\xi_s,\d W_s\>-\ff 1 2 \int_0^t |\xi_s|^2\d s},\ \ t\in [0,T]$$ is a martingale satisfying
\beq\label{NN10} \E[R_t^2|\F_0]\le \e^{c_3 + c_3 \int_0^t\|\mu_s\|_k^2\d s},\ \ t\in [0,T] \end{equation} for some constant $c_3>0$.
By Girsanov's theorem
$$\tt W_t:= W_t-\int_0^t \xi_s\d s,\ \ t\in [0,T]$$ is $m$-dimensional Brownian motion under the probability measure $\Q_T:=R_T\P$.
Since $b^\mu= b^{(1)}+b^{\mu,0}$, we may reformulate the SDE for $\bar X_t$ as
$$ \d \bar X_{t}= b^\mu_t (\bar X_t)\d t + \sigma_t(\bar X_{t}, \Phi^\gamma_{t}\mu)\d \tt W_t,\ \ \bar X_0=X_0, t\in[0,T],$$
so that the weak uniqueness of \eqref{ENN} yields
$\L_{\bar X |\Q_T}=\L_{X^\mu}$. Combining this with \eqref{NN9} and \eqref{NN10}, we obtain
\beg{align*} &\E\big[\e^{\int_0^t f(s,X_s^\mu)^2  \d s}\big|\F_0\big]= \E\big[R_t \e^{\int_0^t f(s,\bar X_s)^2  \d s}\big|\F_0\big] \\
&\le \big(\E[|R_t|^2|\F_0]\big)^{\ff 1 2} \big(\E[ \e^{\int_0^t 2f(s,\bar X_s)^2  \d s}|\F_0]\big)^{\ff 1 2}
 \le \e^{c_4 + c_4\int_0^t\|\mu_s\|_k^2\d s+  c_4 \|f\|_{\tt L^{q}_p(t)}^{c_1}} \end{align*}
 for some constant $c_4>0$. This implies \eqref{JJ1} for some constant $c_0>1$.

  By choosing large enough constant $C_j>0$ such that
 $h(r):= \{\log (C_j+ r) \}^j$ is concave for $r\ge 0$, using Jensen's inequality and \eqref{JJ1} we find a constant $\tt C_j>1$ increasing in $j\ge 1$ such that
\beg{align*} &\E\bigg[\bigg(\int_0^t |f_s(X_s^\mu)|^2\d s\bigg)^j \bigg|\F_0\bigg] \le \E\bigg(\Big[\log \big(C_j+\e^{\int_0^t f_s(X_s^\mu)^2\d s}\big)\Big]^j\bigg|\F_0\bigg) \\
&\le \Big[\log \big(C_j+\E[\e^{\int_0^t f_s(X_s^\mu)^2\d s}]\big|\F_0\big)\Big]^j \le \tt C_j\bigg(1+\int_0^t \|\mu_s\|_k^2\d s+\|f\|_{\tt L^q_p(t)}^{c_1}\bigg)^j.\end{align*}
Using $\ff {f}{\|f\|_{\tt L^q_p(t)}}$ replacing $f$,  we derive
$$ \E\bigg[\bigg(\int_0^t |f_s(X_s^\mu)|^2\d s\bigg)^j\bigg|\F_0\bigg]\le  \|f\|_{\tt L^q_p(t)} ^{2j}   \tt C_j\bigg(1+\int_0^t \|\mu_s\|_k^2\d s+1\bigg)^j$$
which implies \eqref{JJ2}.
\end{proof}

We are now ready to prove Theorem \ref{T1}.

\beg{proof}[Proof of Theorem \ref{T1}]
(1) Since  \eqref{NES} is included in Lemma \ref{LN2}, it
  remains to prove  that  $\Phi^\gg$ has a unique fixed point in $\C_k^{\gg, N}$ for $N>N_0$.

   Under $(A_1)$, \eqref{NN*1} holds for $\kk=0$, so that
      \eqref{WDI} becomes
  $$ \W_{k,\theta} (\Phi^\gg\mu^1,\Phi^\gg\mu^2)
      \le   \bb(\theta)  \W_{k,\theta}(\mu^1,\mu^2),\ \  \theta\ge \theta_0.$$
     Taking large enough  {  $\theta $} such that  $\bb(\theta)<1$    we  prove the contraction of $\Phi^\gg$ on the complete metric space $(\C_k^{\gg,N},\W_{k,\theta})$,
 so that $\Phi^\gg$ has a unique fixed point in $\C_k^{\gg,N}$.

(2) Let $\kk=0$ and $N>0$. For any two solutions  $X_t^i$  of  $\eqref{E1}$ with $\E[|X_0^i|^k]\le N$,
    they   solve \eqref{ENN} for
     $\mu_t^i=\nu_t^i=\L_{X_t^i}, i=1,2$. By \eqref{NES}, there exists a constant $K_N>0$ depending on $N$ such that $\mu,\nu\in \C_k^{\gg,K_N}$.
 Since $\kk=0$ and \eqref{W} for large $\theta$ such that $\ss{c_3\vv_\theta} \le \ff 1 4$, where $\theta$ and $c_3$ depend  on $N$, we obtain
     $$\W_{k,\theta} (\mu_t^1,\mu_t^2)^k\le 2 \ss{c_3} \E[|X_0^1-X_0^2|^k].$$
     Substituting into \eqref{*F} for $\kk=0$ yields the estimate \eqref{NWS} for some constant $C_{j,N}>0.$  When $K_0=0$ we have $|b^{\mu,0}|\le f_0$ for any $\mu\in C([0,T];\scr P_k)$, so that all the above constants  are uniformly bounded in $N$,
     hence  \eqref{NWS} holds for  some constant $C_{j,N}=C_j$ independent of $N$.

     Finally, by taking $j=k$ and $X_0^1, X_0^2$ such that
     $$\L_{X_0^1}=\mu^1,\ \ \L_{X_0^2}=\mu^2,\ \ \E[|X_0^1-X_0^2|^k]=\W_k(\mu^1,\mu^2)^k,$$
     we deduce \eqref{NWS'} from \eqref{NWS}. \end{proof}

\section{Proof of Theorem \ref{T2}  }

By Lemma \ref{LN1}, \eqref{ENN} is well-posed so that the map $\Phi^\gg$ is well-defined on $\C_k^\gg$. Moreover, Lemma \ref{LN2} ensures that
$\C_{k}^{\gg,N}$ is $\Phi^\gg$-invariant for $N\ge N_0$. So, for the well-posedness of \eqref{E1}, it suffices to prove the contraction of $\Phi^\gg$ in $\C_k^{\gg,N}$ for $N>N_0$
under the metric $\W_{k, \kk var,\theta}$ for large $\theta>0$. To this end, we will make use of the parametrix expansion for   transition densities.

 \subsection{Parametrix expansion}

 For any $\mu\in \C_k^\gg$, and a measurable map $\Gamma$ on $\C_k^\gg$,  consider the following
 SDE:
 \beq\label{Et}  \d  {X}_{t}^{x,\mu}= b_t( {X}_{t}^{x,  \mu}, \mu_t)\d t+\si_t( {X}_{t}^{x,\mu},\Gamma_{t}\mu)\d W_{t},\ \ t\in [0,T],\ \ X_{0}^{x,\mu}=x. \end{equation}
Again by \cite[Theorem 2.1]{21Ren},   $(A_0)$ implies the well-posedness of this SDE. Moreover,  by Theorem 6.2.7(ii)-(iii) in \cite{BKRS},    $\L_{X_{t}^{x,\mu}}$ has  a density function
 $p_{t}^{\mu}(x,\cdot)$ (called transition density) with respect to the Lebesgue measure.
 By the standard Markov property of solutions to \eqref{Et}, the solution to \eqref{ENN} satisfies
\beq\label{ES1} \E f(X_{t}^{\mu})= \int_{\R^d}\gg(\d x)  \int_{\R^d}f(y)p_{t}^{\mu}(x,y)\d y,\ \ t\in (0,T], f\in \B_b(\R^d), \end{equation}
where $\B_b(\R^d)$ is the class of bounded measurable functions on $\R^d$. So, to estimate $\|\Phi^\gg_t\mu-\Phi^\gg_t\nu\|_{k,var}$, it suffices to calculate
$|p_t^\mu(x,y)-p_t^\nu(x,y)|$, for which we make use of the parametrix expansion formula.

For  any $x,z\in\R^d, 0\le s<t\le T$ and $\mu\in \C_k^\gg$, let $p_{s,t}^{\mu,z}(x,\cdot)$ be the distribution density function of the random variable
 $$X_{s,t}^{x,\mu,z}:= x + \int_s^t \si_r(z, \Gamma_{r}\mu)\d W_r.$$
  Let
\beq\label{MA}
  a_{s,t}^{\mu,z} :=  \int_s^t (\si_r\si_r^*)(z, \Gamma_{r}\mu)\d r,\ \ 0\le s<t\le T. \end{equation}
We have
\beq\label{ES6'} p_{s,t}^{\mu,z}(x,y)= \ff{\exp[-\ff 1 2 \<(a_{s,t}^{\mu,z})^{-1}(y-x),  y-x\>]}{(2\pi)^{\ff d 2} ({\rm det} \{a_{s,t}^{\mu,z}\})^{\ff 1 2}},\ \ x,y\in \R^d.\end{equation}
Obviously, $(A_0)$ and $(A_2)$   imply
\beq\label{mv}\beg{split}
 &\|a_{s,t}^{\mu,z}-a_{s,t}^{\nu,z}\|\le K\int_s^t  \W_{k}\big(\Gamma_{r}\mu, \Gamma_{r}\nu\big)\d r,\\
 &\ff{1}{K(t-s)}\le \|(a_{s,t}^{\mu,z})^{-1}\|\leq \ff{K}{t-s},\ \ 0\le s<t\le T,\ \mu,\nu\in \C_k^\gg.\end{split}
\end{equation}
Next, for $\mu\in \C_k^\gg$, $y,z\in\R^d$ and $0\leq s<t\leq T$, let
\beq\label{ES7}\beg{split} &H_{s,t}^{\mu,1}(y,z)= H_{s,t}^{\mu} (y,z):=    \left\<-b_{s}(y, \mu_s),\nabla p_{s,t}^{\mu,z}(\cdot,z)(y)\right\>\\
&\qquad +\frac{1}{2}\mathrm{tr}\left[\left\{(\sigma_{s}\sigma^\ast_{s})(z, \Gamma_{s}\mu)-(\sigma_{s}\sigma^\ast_{s})(y, \Gamma_{s}\mu)\right\}\nabla^2 p_{s,t}^{\mu,z} (\cdot,z)(y)\right],\\
&H_{s,t}^{\mu,j} (y,z) := \int_s^t\d r \int_{\R^d}  H_{r,t}^{\mu, j-1}(z',z) H_{s,r}^{ \mu}(y,z') \d z',\ j\ge 2.\end{split}\end{equation}
By the parabolic equations for the transition densities $p_{s,t}^{\mu}$ and $p_{s,t}^{\mu,z}$, see   for instance the paragraph after  Lemma 3.1 in \cite{KM00},  we have the parametrix expansion formula
\beq\label{ES8} p_{t}^{\mu}(x,z) = p_{0,t}^{\mu,z}(x,z) + \sum_{j=1}^\infty \int_0^t \d s \int_{\R^d} H_{s,t}^{\mu,j}  (y,z) p_{0,s}^{\mu,z} (x,y) \d y.\end{equation}

Let
 \beq\label{ES10} \tt p_{s,t}^K(x,y)= \ff{\exp[-\ff 1 {4K(t-s)} |y-x|^2]}{(4K\pi (t-s))^{\ff d 2}},\ \ x,y\in \R^d, 0\le s<t\le T.\end{equation}
By multiplying the time parameter with $T^{-1}$ to make it stay in $[0,1]$, we deduce from \cite[(2.3), (2.4)]{ZG} with $\beta=\beta'=1$ and $\ll= \ff 1 {8KT}$  that
\begin{equation}\label{CK}\begin{split} &\int_s^t\int_{\R^d} \tt p^K_{s,r}(x,y')(r-s)^{-\frac{1}{2}} g_r(y')(t-r)^{-\frac{1}{2}} \tt p^{2K}_{r,t}(y',y) \d y'\\
&\leq c  (t-s)^{-\ff 1 2 +\frac{1}{2}(1-\frac{d}{p}-\frac{2}{q})}\tt p^{2K}_{s,t}(x,y) \|g\|_{\tilde{L}_p^q([s,t])},\ \ 0\leq s<t\leq T, g\in \tt L_p^q([s,t])\end{split}\end{equation}
holds for some constant $c>0$ depending on $T,d,p,q$ and $K$.
By the condition on $a$ included in $(A_0)$,    we find a constant $c_1>0$ such that  \eqref{ES6'} implies
\beq\label{G2}\beg{split} &p_{s,t}^{\mu,z}(x,y) \bigg(1+ \ff{|x-y|^4}{(t-s)^2} \bigg)\\
& \le c_1   \tt p^K_{s,t}(x,y),\ \ x, y,z\in\R^d, 0\le s< t\le T, \gg\in \scr P_k,\mu\in C([s,t];\scr P_k).\end{split}\end{equation}

\beg{lem}\label{L2}  Assume $(A_0)$ and $(A_2)$.  Let $p_{s,t}^{\mu,z}(x,y)$ be defined by $ \eqref{MA}$ and  $\eqref{ES6'}$ for some map $\GG: \C_k^\gg\to\C_k^\gg$.  There exists a constant $c>0$ independent of $\GG$, such that for any $0\le s<t\le T, x,y,z\in \R^d, \gg\in\scr P_k,$ and $
\mu,\nu\in C([s,t];\scr P_k)$,
\begin{align}\label{g0}
&\bigg(1+\ff{|x-y|^2}{t-s}\bigg) |p_{s,t}^{\mu,z}(x,y)-p_{s,t}^{\nu,z}(x,y)|\leq \ff{c \tt p^K_{s,t}(x,y)}{t-s}\int_s^t \W_{k}(\Gamma_{r}\mu, \Gamma_{r}\nu)\d r,
\end{align}
\beq\label{g1} \ss{t-s} |\nabla p_{s,t}^{\mu,z}(\cdot,y)(x)|+ (t-s)  \|\nabla^2 p_{s,t}^{\mu,z}(\cdot,y)(x)\|\leq c\tt p^K_{s,t}(x,y),\end{equation}
\beq\label{g2}\beg{split}
&\ss{t-s} |\nabla p_{s,t}^{\mu,z}(\cdot,y)(x)-\nabla p_{s,t}^{\nu,z}(\cdot,y)(x)|\\
&+(t-s) \|\nabla ^2p_{s,t}^{\mu,z}(\cdot,y)(x)-\nabla^2 p_{s,t}^{\nu,z}(\cdot,y)(x)\|\\
 &\leq   \ff{c\tt p^K_{s,t} (x,y)}{t-s}  \int_s^t \W_{k}(\Gamma_{r}\mu, \Gamma_{r}\nu)\d r.
\end{split}\end{equation}
\end{lem}

\beg{proof} (1) For fixed $x,y\in\R^d$ and $0\le s<t\le T$, let
$$F(\mu):=\<(a_{s,t}^{\mu,z})^{-1}(y-x),  y-x\>,\ \ \mu\in C([s,t];\scr P_k).$$  It is easy to  see that
\beq\label{LL}\begin{split}
&|p_{s,t}^{\mu,z}(x,y)- p_{s,t}^{\nu,z}(x,y)|\\
&=\left| \ff{\exp[-\ff 1 2 F(\mu)]}{(2\pi)^{\ff d 2} ({\rm det} \{a_{s,t}^{\mu,z}\})^{\ff 1 2}} -\ff{\exp[-\ff 1 2 F(\nu)]}{(2\pi)^{\ff d 2} ({\rm det} \{a_{s,t}^{\nu,z}\})^{\ff 1 2}}\right|\le I_1+I_2,\end{split}\end{equation} where
\beg{align*}& I_1:=   \ff{\left|\exp[-\ff 1 2 F(\mu)]-\exp[-\ff 1 2 F(\nu)]\right|}{(2\pi)^{\ff d 2} ({\rm det} \{a_{s,t}^{\mu,z}\})^{\ff 1 2}}\\
&I_2:= \ff{\exp[-\ff 1 2 F(\nu)]}{(2\pi)^{\ff d 2}}\left| ({\rm det} \{a_{s,t}^{\mu,z}\})^{-\ff 1 2}-({\rm det} \{a_{s,t}^{\nu,z}\})^{-\ff 1 2}\right|.
\end{align*}
 Combining this with  $(A_0)$ and $(A_2)$ which imply   \eqref{mv}, we find a constant $c_1>0$ such that
\begin{align*}
&\left|F(\mu)-F(\nu)\right|= \big| \<\{(a_{s,t}^{\mu,z})^{-1}- (a_{s,t}^{\nu,z})^{-1}\} (y-x),  y-x\>\big|  \\
 &\le c_1   \ff{|y-x|^2}{(t-s)^2}    \int_s^t \W_{k}(\Gamma_{r}\mu, \Gamma_{r}\nu)\d r,
\end{align*}
which together with \eqref{G2} and $\ff{|x-y|^2}{t-s}\leq \frac{1}{2}(1+\ff{|x-y|^4}{(t-s)^2})$ yields that for some constant $c_2>0$,
\begin{equation*}\begin{split}\bigg(1+\ff{|x-y|^2}{t-s}\bigg)I_1
& \leq \ff{c_2  \tt p^K_{s,t}(x,y) }{t-s}  \int_s^t \W_{k}(\Gamma_{r}\mu, \Gamma_{r}\nu)\d r.\end{split}\end{equation*}
Again by \eqref{mv}, \eqref{G2} and $\ff{|x-y|^2}{t-s}\leq \frac{1}{2}(1+\ff{|x-y|^4}{(t-s)^2})$,  we find a constant $c_3>0$ such that
\begin{equation*} \bigg(1+\ff{|x-y|^2}{t-s}\bigg)I_2
 \leq \ff{c_3  \tt p^K_{s,t}(x,y) }{t-s}  \int_s^t \W_{k}(\Gamma_{r}\mu, \Gamma_{r}\nu)\d r.
 \end{equation*}
Combining these with \eqref{LL}, we arrive at
\beg{align*} &\bigg(1+\ff{|x-y|^2}{t-s}\bigg)|p_{s,t}^{\mu,z}(x,y)- p_{s,t}^{\nu,z}(x,y)|\le \ff{(c_2+c_3)  \tt p^K_{s,t}(x,y)}{t-s}  \int_s^t \W_{k}(\Gamma_{r}\mu, \Gamma_{r}\nu)\d r.\end{align*}

(2) By \eqref{ES6'} we have
\beq\label{GP1}
\nabla p_{s,t}^{\mu,z}(\cdot,y)(x)=(a_{s,t}^{\mu,z})^{-1}(y-x)p_{s,t}^{\mu,z}(x,y),\end{equation}
\beq\label{GP2}\beg{split}
\nabla^2p_{s,t}^{\mu,z}(\cdot,y)(x)&=p_{s,t}^{\mu,z}(x,y) \Big( \big\{(a_{s,t}^{\mu,z})^{-1}(y-x)\big\} \otimes\big\{(a_{s,t}^{\mu,z})^{-1}(y-x)\big\}
-(a_{s,t}^{\mu,z})^{-1}\Big)\end{split}. \end{equation}
So, by \eqref{mv} and \eqref{G2} we find a constant $c>0$ such that \eqref{g1} holds.
Moreover,   \eqref{GP1} implies
\begin{align*}
&|\nabla p_{s,t}^{\mu,z}(\cdot,y)(x)-\nabla p_{s,t}^{\nu,z}(\cdot,y)(x)|\\
&\leq \left|\{(a_{s,t}^{\mu,z})^{-1}-(a_{s,t}^{\nu,z})^{-1}\} (y-x)\right|p_{s,t}^{\mu,z}(x,y)\\
&\quad +\big| p_{s,t}^{\mu,z}(x,y)-p_{s,t}^{\nu,z}(x,y) \big|\cdot \big|(a_{s,t}^{\nu,z})^{-1}(y-x)\big|.\end{align*}
Combining this with \eqref{mv}, \eqref{G2} and \eqref{g0}, we find a constant $c>0$ such that
$$|\nabla p_{s,t}^{\mu,z}(\cdot,y)(x)-\nabla p_{s,t}^{\nu,z}(\cdot,y)(x)|\le \ff {c   \tt p^K_{s,t}(x,y)}{(t-s)^{\ff 3 2}}  \int_s^t \W_{k}(\Gamma_{r}\mu, \Gamma_{r}\nu)\d r.$$
Similarly, combining  \eqref{GP2} with  \eqref{mv}, \eqref{G2} and \eqref{g0},   we find a constant $c>0$ such that
$$\|\nabla^2 p_{s,t}^{\mu,z}(\cdot,y)(x)-\nabla^2 p_{s,t}^{\nu,z}(\cdot,y)(x)\|\le \ff {c   \tt p^K_{s,t}(x,y)}{(t-s)^2}  \int_s^t \W_{k}(\Gamma_{r}\mu, \Gamma_{r}\nu)\d r.$$
 Therefore, \eqref{g2} holds for some constant $c>0$.
\end{proof}

For $ 0\le s\le  t\le T,\gg\in\scr P_k $ and $\mu,\nu\in C([s,t];\scr P_k)$, let
 \beq\label{ga} \LL_{s,t}(\mu,\nu):=\sup_{r\in[s,t]}\big\{ \W_{k}(\Gamma_{r}\mu, \Gamma_{r}\nu) +\W_{k,\kk var}(\mu_r,\nu_r)\big\}. \end{equation}
\beg{lem}\label{L3} Assume  $(A_0)$ and $(A_2)$.  Let $\dd:= \frac{1}{2}\left(1-  \frac{d}{p_0} -   \frac{2}{q_0} \right)>0$ and denote
$$
 S_\mu:= \sup_{t\in \in [0,T]} (1+\|\mu_t\|_k),\ \
 S_{\mu,\nu}:=S_{\mu}\lor S_{\nu},\ \ \nu,\mu \in \C_k^\gg.$$
 Then there exists a constant
  $C\ge 1$ such that for any $0\le s<t\le T$,  $ y,z\in\R^d$,    $\mu,\nu\in \C_k^\gg$, and $j\ge 1$,
\beq
\label{Ta2'}  | H_{s,t}^{\mu,j}(y,z) |  \le  f_0(s,y) (CS_\mu)^j (t-s)^{-\frac{1}{2}+\delta (j-1)} \tt p^{2K}_{s,t}(x,y),
\end{equation}
\beq
 \label{Ta3'} \beg{split} &|H_{s,t}^{\mu,j}(y,z) -H_{s,t}^{\nu,j}(y,z) |\\
 &\le  j f_0(s,y) (C S_{\mu,\nu})^j (t-s)^{-\frac{1}{2}+\delta (j-1)}\tt p^{2K}_{s,t}(x,y) \LL_{s,t}(\mu,\nu).
\end{split} \end{equation}
\end{lem}

\beg{proof}   (1) By \eqref{ES7}, \eqref{g1}, $(A_0)$ and $(A_2)$,    we  find  a constant $c_1>0$   such that for any  $0\le s<t\le T, \mu \in C([0,T];\scr P_k)$ and $y,z\in \R^d$,
\beq\label{A1} \begin{split}
 |H_{s,t}^{\mu}(y,z)|\le c_1(t-s)^{-\frac{1}{2}}\{(1+\|\mu_s\|_k)f_0(s,y)\} \tt p^{K}_{s,t}(y,z).  \end{split} \end{equation}
So, \eqref{Ta2'} holds for $j=1$ and $C=c_1$.
Thanks to \cite[(2.3), (2.4)]{ZG} with $\beta=\beta'=1$, $\lambda=\frac{1}{8K}$, we have
\begin{align}\label{0*2}  \nonumber I_j :=&\int_s^t\int_{\R^d} (t-u)^{-\frac{1}{2}}(t-u)^{\delta(j-1)}\tt p^{2K}_{u,t}(y,z)f_0(u,y)(u-s)^{-\frac{1}{2}} \tt p^{K}_{s,u}(x,y)\d y\d u\\
 &\leq c_2(t-s)^{-\frac{1}{2}}\tt p^{2K}_{s,t}(x,z)(t-s)^{\frac{1}{2}(1-\frac{d}{p_0}-\frac{2}{q_0})}\|f_0\|_{\tilde{L}_{p_0}^{q_0}([s,t])}(t-s)^{\delta(j-1)}\\
 \nonumber&=c_3(t-s)^{-\frac{1}{2}}\tt p^{2K}_{s,t}(x,z)(t-s)^{\delta j}.  \ \ 0\le s<t\le T,j\geq 1
\end{align}
where $c_3:=   c_2\|f_0\|_{\tilde{L}_p^q([s,t])}$.
   Let  $C:= 1\lor c_1^2\lor (4c_3^2)$.  If for some $j\ge 1$ we have
\begin{align*}
| H_{s,t}^{\mu,j}(y,z) |   &\le (CS_\mu)^j  f_0(s,y)\tilde{p}^{2K}_{s,t}(y,z)
  (t-s)^{-\frac{1}{2}+\delta(j-1)}
\end{align*} for all $y,z\in\R^d$ and $0\le s<t\le T$,
then by combining   with \eqref{A1} and \eqref{0*2}, we arrive at
\beg{align*}  | H_{s,t}^{\mu,j+1}(y,z) | &\le  \int_s^t\d u \int_{\R^d}  |H_{u,t}^{\mu,j}(z',z) H_{s,u}^{\mu}(y,z') | \d z'\\
&\le C^j \ss C (S_\mu)^{j+1}  f_0(s,y)I_k\\
&\le  C^{j+1}(S_\mu)^{j+1}  f_0(s,y)(t-s)^{-\ff 1 2 +\dd j} \tilde{p}^{2K}_{s,t}(y,z) .
\end{align*}
Therefore, \eqref{Ta2'} holds for all $j\ge 1$.

(2) By    \eqref{g1}, \eqref{g2}, \eqref{mv}, $(A_0)$ and $(A_2)$, we find a constant $c>0$ such that  for any  $0\le s<t\le T, \mu,\nu\in C([0,T];\scr P_k)$ and $y,z\in \R^d$,
\beq\label{A2}| H_{s,t}^{\mu} (y,z)- H_{s,t}^{\nu} (y,z)|
  \leq  c(t-s)^{-\ff 1 2}\tt p^{K}_{s,t}(y,z)S_{\mu,\nu}f_0(s,y) \LL_{s,t}(\mu,\nu).
 \end{equation}
Let, for instance, $L= 1+ 4C^2 + 4 c^2$, where $C$ is in \eqref{Ta2'}. If for some $j\ge 1$ we have
\begin{align*}&| H_{s,t}^{\mu,j}(z',z)-H_{s,t}^{\nu,j}(z',z)| \le j(LS_{\mu,\nu})^j f_0(s,z')  \tilde{p}^{2K}_{s,t}(z',z)(t-s)^{-\frac{1}{2}+\delta(j-1)}\LL_{s,t}(\mu,\nu),
\end{align*}  for any $0\le s<t\le T$ and $z,z'\in\R^d$, then \eqref{Ta2'}, \eqref{0*2} and  \eqref{A2} imply
\beg{align*} & | H_{s,t}^{\mu,j+1}(y,z) - H_{s,t}^{\nu,j+1}(y,z)|  \\
& \le   \int_s^t\d r \int_{\R^d} \Big\{|H_{r,t}^{\mu, j}(z',z)-H_{r,t}^{\nu, j}(z',z)|\cdot | H_{s,r}^{\mu}(y,z') |  \\
&\qquad\qquad  \qquad \quad+ |H_{r,t}^{\nu, j}(z',z)  | \cdot| H_{s,r}^{\mu}(y,z')-H_{s,r}^{\nu}(y,z')|\Big\}\d z'\\
&\leq  (j+1)(LS_{\mu,\nu})^{j+1}f_0(s,y)  \tilde{p}^{2K}_{s,t}(y,z)(t-s)^{-\frac{1}{2}+\delta j}\LL_{s,t}(\mu,\nu).\end{align*}
Therefore,   \eqref{Ta3'} holds for  some constant $C>0. $
\end{proof}

 We are now ready to prove  the following main result in this part, which will be used to prove the contraction of $\Phi^\gg$
on the path space over a small time interval. For $t_0\in (0,T]$, let
$$\C_{k,t_0}^{\gg,N}:= \big\{\mu\in C([0,t_0];\scr P_k):\ \mu_{\cdot\land t_0} \in \C_k^{\gg, N}\big\},\ \ N\ge N_0.$$

 \beg{lem}\label{LM} Assume  $(A_0)$ and $(A_2)$.   For any $N\ge N_0$, there exist $\theta_N>0, t_N\in (0,T]$   such that
$$\W_{k,\kk var,\theta_N }(\Phi_{\cdot\land t_N}^\gg\mu,  \Phi_{\cdot\land t_N}^\gg\nu)  \le  \ff 1 2  \W_{k,\kk var,\theta_N}(\mu_{\cdot\wedge t_N},\nu_{\cdot\wedge t_N}),\ \ \mu,\nu\in \C_{k,t_N}^{\gg,N}.$$
 \end{lem}

\beg{proof}
 By \eqref{G2}, Lemma \ref{L2}, Lemma \ref{L3},  \eqref{ES8}, \eqref{CK} and $(A_2)$,    we find  constants $c_1,c_2,c_3>0$  such that for any $\theta>0$ and $ t_N\in (0, T\land (2CN)^{-\ff 1 \dd}],$
  \beg{align*}& |p_t^{\mu}(x,z)-p_{t}^{\nu}(x,z)|
   \leq \ff{c_1 \tt p^K_{0,t}(x,z)}t \int_0^t \W_{k}(\Gamma_{s}\mu, \Gamma_{s}\nu)\d s\\
   &\qquad + \sum_{n=1}^\infty \int_0^t  \d s \int_{\R^d}\big\{ |H_{s,t}^{\mu,n}   -H_{s,t}^{\nu,n} | (y,z)   p_{0,s}^{\nu,z} (x,y)+
      |H_{s,t}^{\mu,n}  (y,z)| |p_{0,s}^{\mu,z}  -p_{0,s}^{\nu,z} |(x,y) \big\}\d y \\
  &\leq  c_1  \e^{\theta  t} \W_{k,\theta}(\Gamma_{\cdot\wedge t}\mu, \Gamma_{\cdot\wedge t}\nu)\tt p^K_{0,t}(x,z) \\
  &+\sum_{n=1}^\infty (n+1) (C N)^n \Lambda_{0,t}(\mu,\nu) t^{\frac{1}{2}+\delta (n-1)}\\
  &\qquad  \times \int_0^t\int_{\R^d}(t-r)^{-\frac{1}{2}}\tt p^{2K}_{r,t}(y,z)f_0(r,y)r^{-\frac{1}{2}}\tt p^{K}_{0,r}(x,y)\d y \d r\\
  &\leq  c_1  \e^{t\theta}   \W_{k,\theta}(\Gamma_{\cdot\wedge t}\mu, \Gamma_{\cdot\wedge t}\nu)\tt p^K_{0,t}(x,z) +c_2t^\dd \Lambda_{0,t}(\mu,\nu) \tt p_{0,t}^{2K}(x,z)\sum_{n=1}^\infty (n+1) (C N)^n t^{ \delta (n-1)}\\
  &\leq  c_1  \e^{\theta t}  \W_{k,\theta}(\Gamma_{\cdot\wedge t}\mu, \Gamma_{\cdot\wedge t}\nu)\tt p^K_{0,t}(x,z)+ c_3 t^\delta \Lambda_{0,t}(\mu,\nu)\tt p_{0,t}^{2K}(x,z)  \end{align*}
  holds for any $ x,z\in\R^d, t\in (0,t_N],  \mu,\nu\in \scr C_k^{\gg,N}.$
Combining this with \eqref{ES10}, we find a constant $c_4>0$ such that
\beq\label{R*S} \begin{split}
& \sup_{|g|\leq 1+|\cdot|^k}\left|\int_{\R^d}\int_{\R^d}g(z)(p_{t}^{\mu}- p_{t}^{\nu}) (x,z)\d z\gamma(\d x) \right|\\
&\le  c_1     \e^{\theta t}  \W_{k,\theta}(\Gamma_{\cdot\wedge t}\mu, \Gamma_{\cdot\wedge t}\nu)\int_{\R^d\times\R^d} (1+|z|^k)
 \tt p_{0,t}^{K}(x,z) \d z\gg(\d x)\\
 &\qquad +c_3t^\delta\Lambda_{0,t}(\mu,\nu)  \int_{\R^d\times\R^d} (1+|z|^k)
 \tt p_{0,t}^{2K}(x,z) \d z\gg(\d x)\\
 &\le c_4 \e^{\theta t}  \W_{k,\theta}(\Gamma_{\cdot\wedge t}\mu, \Gamma_{\cdot\wedge t}\nu)+  c_4t^\delta\Lambda_{0,t}(\mu,\nu) ,\ \  t\in (0,t_N],\ \mu,\nu\in \C_k^{\gg,N}.
 \end{split}\end{equation}
 Taking $\Gamma=\Phi^\gamma$,  by the definition  of $\Phi_{t}^\gg$,  \eqref{R*S}  and \eqref{ga}, we find a constant $c_5>0$ such that
 \begin{align*}  &\W_{k,\kappa var,\theta}(\Phi_{\cdot\land t_N}^\gg\mu, \Phi_{\cdot\land t_N}^\gg\nu)\\
&\leq c_5\W_{k,\theta}(\Phi_{\cdot\wedge t_N}^\gg\mu, \Phi_{\cdot\wedge t_N}^\gg\nu)+  c_5t_N^\delta\W_{k,\kk var,\theta}(\mu_{\cdot\wedge t_N},\nu_{\cdot\wedge t_N}),\ \  \ \mu,\nu\in \C_k^{\gg,N}.
\end{align*}
By \eqref{WDI} with $\bb(\theta)\to 0$ as $\theta\to\infty$, we find large $\theta_N>0$ and small $t_N\in (0,T]$ depending  on $N$ such that
 \begin{align*}\W_{k,\kappa var, \theta_N}(\Phi_{\cdot\land t_N}^\gg\mu, \Phi_{\cdot\land t_N}^\gg\nu)&\le  c_5(\bb(\theta_N) +t_N^\delta) \W_{k,\kk var,\theta_N}(\mu_{\cdot\wedge t_N},\nu_{\cdot\wedge t_N})\\
&\le \ff 1 2\W_{k,\kk var,\theta_N}(\mu_{\cdot\wedge t_N},\nu_{\cdot\wedge t_N}).
\end{align*}
\end{proof}

\subsection{Proof of Theorem \ref{T2} }

Estimate \eqref{NES} is included in Lemma \ref{LN2}(2). It suffices to prove the well-posedness of \eqref{E1} and estimate \eqref{NES2} for $\kk=0,$ where $C(N)$ is bounded in $N$
when $K_0=0$.

(a) Well-posedness.
By the priori estimate  \eqref{NES}, there exists a constant $C>0$ such that for any solution of \eqref{E1}
 on $[0,T]$ with $\L_{X_0}=\gg$,
\beq\label{A*B}
  \sup_{t\in [0,T]}\L_{X_t}(|\cdot|^k)\le C.\end{equation}  So, we may fix $N_0>0$ depending only on $C$ such that  any solution of \eqref{E1} with initial distribution $\gg$ satisfies
$\L_{X_\cdot}\in \C_k^{\gg, N_0}.$  By Lemma \ref{LM},   there exists $\theta>0$ and $t_0\in (0,T]$ depending only on $N_0$ such that the map  $\Phi^\gamma_{\cdot\land t_0}$ is contractive in $\C_{k,t_0}^{\gg, N_0}$
under the metric $\W_{k, \kk var,\theta}$, and hence   \eqref{E1} for $t\in [0,t_0]$ is well-posed for distributions in $\scr P_k$ and   \eqref{A*B} holds.
Using $(t_0,X_{t_0})$  replacing $(0, X_0)$, the same argument implies the well-posedness of \eqref{E1} for $t\in [t_0, (2 t_0)\land T]$ and that \eqref{A*B} holds for $(2 t_0)\land T$ replacing $t_0$.
 By repeating the procedure finitely many times, we prove the well-posedness of \eqref{E1} for distributions
 in $\scr P_k$.

  (b) Estimate \eqref{NES2}.
  For any  $\mu^i_0 \in\scr P_{k}$ with $\ \mu^i_0(|\cdot|^k)\leq N, i=1,2,$ let
    $$\mu_t^i=P_t^\ast\mu_0^i, \ \ i=1,2, t\in [0,T].$$  By \eqref{NES}, there exists a constant $C_N>0$ such that
  \beq\label{RS0} \sup_{t\in [0,T]}(\mu_t^1+\mu_t^2)(|\cdot|^k)\le C_N.\end{equation}
 So, there exists a constant $\bar{N}$ depending on $C_N$ such that $$\mu_\cdot^i\in \scr C _{k}^{\gamma,\bar{N}}, \ \ i=1,2.$$
Consider the SDEs
\begin{align}\label{AEQ}
\d X_t^{x,i}=b_t(X_t^{x,i},\mu^i_t)+\sigma_t(X_t^{x,i},\mu_t^i)\d W_t,\ \ X_0^{x,i}=x\in\R^d, t\in [0,T], i=1,2.
\end{align}
We have
\beq\label{RS} \mu_t^i:=P_t^*\mu^i_0= \int_{\R^d} \L_{X_t^{x,i}}\mu_0^i(\d x),\ \ t\in [0,T], i=1,2,\end{equation}
According to   \cite[Theorem 2.1(2)]{W21d},   \eqref{RS0} and $(A_0)$ imply
$$\|\L_{X_t^{x,i}}-\L_{X_t^{y,i}}\|_{var} \le \ff{c_1}{\ss t} |x-y|,\ \ x,y\in \R^d, t\in (0,T], i=1,2$$ for some constant $c_1>0$ depending on $N$.
Combining this with \eqref{RS} gives
  \begin{align}\label{var}\nonumber&\bigg\|P_t^*\mu_0^1- \int_{\R^d}\L_{X_t^{y,1}}\mu_0^2(\d y)\bigg\|_{var}= \left\|\int_{\R^d}\L_{X_t^{x,1}}\mu_0^1(\d x)-\int_{\R^d}\L_{X_t^{y,1}}\mu_0^2(\d y)\right\|_{var}\\
&\leq \inf_{\pi\in\C(\mu_0^1,\mu_0^2)}\int_{\R^d\times \R^d}\|\L_{X_t^{x,1}}- \L_{X_t^{y,1}}\|_{var} \pi(\d x,\d y)\\
\nonumber&\leq  \frac{c_1}{\sqrt{t}}\W_1(\mu_0^1,\mu_0^2)\leq  \frac{c_1}{\sqrt{t}}\W_k(\mu_0^1,\mu_0^2),\ \ t\in (0,T].
\end{align}
 On the other hand, by \eqref{RS} and  \eqref{R*S} for $\mu=\mu^1$, $\nu=\mu^2$, $\kk=0$ and $\Gamma=id$, we find   constants $c_2>0$ and $t_N\in (0,T]$ depending on $N$ such that
\begin{align*}&\left\|P_t^\ast\mu_0^2-\int_{\R^d}\L_{X_t^{y,1}}\mu_0^2(\d y)\right\|_{var}
  \leq c_2\sup_{t\in [0,T]}\W_k(\mu_t^1,\mu_t^2),\ \ t\in [0,t_N].
\end{align*}
  For any $t\in [t_N, T]$, repeating the above argument for the time interval $[t-t_N, t]$ replacing $[0,t_N]$  we prove
\begin{align*}&\left\|P_t^\ast\mu_0^2-\int_{\R^d}\L_{X_t^{y,1}}\mu_0^2(\d y)\right\|_{var}
  \leq c\sup_{t\in [0,T]}\W_k(\mu_t^1,\mu_t^2)
\end{align*}
for some constant $c>0$ depending on $N$.
  Combining this with \eqref{var} and \eqref{NWS'} which holds since $(A_2)$ with $\kk=0$ implies $(A_1),$ we prove \eqref{NES2} for some constant $C(N)>0.$

Finally, noting that the dependence on $N$ comes from Krylov's and Khasminskii's estimates for
the solutions, and when $K_0=0$ we have $|b^{\mu,0}|\le f_0$   for all $\mu\in \C_k$,  these estimates are uniform in $\mu$.
Thus, in this case \eqref{NES2} holds for all $\mu,\nu\in \scr P_k$ and a constant $C>0$ independent of $N$.

\section{Extension of Theorem \ref{T1} to reflecting SDEs}

Let $D\subset \R^d$ be a connected  open domain with   $\pp D\in C_b^{2,L}$ in the following sense:
there exists a constant $r_0>0$ such that  the polar coordinate map
$$ \Psi: \pp D\times [-r_0,r_0] \ni (z, r)\mapsto   z+ r\n(z) \in \pp_{\pm r_0}D:=\big\{x\in\R^d: \rr_\pp(x):={\rm dist}(x,\pp D)\le r_0\big\} $$
is a $C^2$-diffeomorphism, such that   $ \Psi^{-1}(x)$  have   bounded and continuous first and second order derivatives in $x\in  \pp_{\pm r_0}D$,   and  $\nn^2\rr_\pp$ is Lipschitz continuous on $\pp_{\pm r_0}D.$

 Consider the following  distribution dependent   reflecting SDE on the closure $\bar D$ of $D$:
\beq\label{E1'} \d X_t= b_t(X_t, \L_{X_t}) \d t+ \si_t(X_t,\L_{X_t})\d W_t + \n(X_t) \d l_t,\ \ t\in [0,T],\end{equation}
where     $\n$ is the  unit inward normal vector field on the boundary $\pp D$ and $l_t$ is a continuous adapted increasing process with $\d l_t$ supported on $\{t: X_t\in\pp D\}.$
 Let $\tt L_p^q(T)$ and $\scr P_k$ be defined as before for $\bar D$ replacing $\R^d$. When $\si_t(x,\mu)=\si_t(x)$ does not depend on $\mu$, the well-posedness of \eqref{E1'} has been proved in \cite{W21a} under the following assumption, where $a_t^\mu:= (\si_t\si_t^*)(\cdot,\mu_t)$.

 \beg{enumerate} \item[$(B)$] Assumptions $(A_0)$ and $(A_1)$    hold  for $\bar D$ replacing $\R^d$. Moreover, there exists a constant $c>0$ such that for any $\mu\in C([0,T];\scr P_k),$ the Neumann semigroup $\{P^\mu_{s,t}\}_{0\le s\le t\le T}$  generated by  the operator
 $L_t^\mu:=\ff 1 2 {\rm tr}\{a_t^\mu\nn^2\}+ b_t^{(1)}\cdot\nn $ on $\bar D$
 satisfies
 \beq\label{AB1}  \|\nn^i P_{s,t}^{\mu} \phi\|_\infty \le  c (t-s)^{-\ff i 2} \|\phi\|_\infty,\ \ 0\le s<t\le T, \phi\in C^i_b(\bar D),\ i=1,2.\end{equation}
\end{enumerate}
\beg{thm} Assume   $(B)$ and let $\pp D\in C_b^{2,L}.$  Then the assertions in Theorem \ref{T1} hold for $\eqref{E1'}$ replacing $\eqref{E1}.$\end{thm}

\beg{proof} Let $\gg\in \scr P_k$ and consider the initial value $X_0$ with $\L_{X_0}=\gg$. It suffices to prove that Lemmas \ref{LN1}-\ref{LJJ} hold for  $\Phi^\gg\mu:=\L_{X^\mu}$  with  the following reflecting SDE replacing \eqref{ENN}:
\beq\label{*ENN}   \d X_t^{\mu}= b_t(X_t^{\mu}, \mu_t)\d t+\si_t(X_t^{\mu},\L_{X_t^{\mu}})\d W_{t}+\n(X_t^\mu) \d l_t^\mu,\ \
t\in [0,T], X_0^{\mu}=X_0.\end{equation}

(a) Assertions in Lemma \ref{LN1}. For $\gg^i\in \scr P_k, \mu^i\in \C_k^{\gg^i,N}$ and $\nu^i\in \C_k^{\gg^i}$, $i=1,2$, instead of \eqref{W*1} we consider the reflecting SDEs
\begin{equation*}  \d X_{t}^{i}= b_t^{\mu^i}(X_t^{i})\d t+\si_t^{\nu^i}(X_t^{i})\d W_{t}+\n(X_t^{i})\d l_t^i,\  \  \L_{X_0^i}=\gg^i, t\in [0,T], i=1,2. \end{equation*}
By \cite[Theorem 2.2(ii)]{W21a}, $(B)$ implies the well-posedness of these reflecting SDEs.

Next, according to the proof of \cite[Theorem 2.2(ii)]{W21a}, there exists a semimartingale $H_t$ such that
$$C^{-1} |X_t^1-X_t^2|^2\le H_t\le C|X_t^1-X_t^2|^2,\ \ t\in [0,T]$$ holds for some constant $C>1$, and  instead of \eqref{**1},
$$  \d H_t^{j}\le c_2  \eta_t^{2j} \d \{A_t + l_t^1+l_t^2\}+ c_2 (g_t^2+\tt g_t) \big\{\W_{k}(\mu^1_t,\mu_t^2)^{2j} + \W_k(\nu^1_t,\nu_t^2)^{2j} \big\}\d t +\d M_t$$
holds for some constant $c_2>0$ and all $t\in [0,T]$.

Then the desired assertions can be proved as in the proof of Lemma \ref{LN1}  by using Khasminskii's estimate in \cite[Lemma 2.7]{W21a}, as well as the  estimate
\beq\label{F*} \E \big[\e^{\ll (l_T^1+l_T^2)}\big]\le \e^{c(1+\ll^2)},\ \ \ll>0\end{equation} for some constant $c>0$  presented in   \cite[Lemma 2.5]{W21a},  where condition $(A_0^{a,b})$ follows from $(A_0)$
included in $(B)$,  according to \cite[Lemma 2.6]{W21a}.

(b) Proof of Lemma \ref{LN2}. In the present case \eqref{**3} becomes
$$    |X_t|^{2j} -|X_0|^{2j}
 \le    c_1 \int_0^t \big\{ 1+ |X_s|^{2j} + |X_s|^{2j-1}f_0(s,X_s)+ \|\mu_s\|_k^{2j}  \big\}\d s +c_1\int_0^t |X_s|^{2j-1}\d l_s+ M_t,$$
 such that \eqref{**4} reduces to
 $$ \E\big(\eta_t\big|\F_0\big)\le c_5 + c_5 |X_0|^{2j}+c_5 \int_0^t \big\{\E(\eta_s|\F_0)+ \|\mu_s\|_k^{2j}\big\}\d s + c_5 \int_0^t\E(\eta_s|\F_0)\d l_s,\ \ t\in [0,T].$$
 Combining this with \eqref{F*} for $l_T^1=l_T$ and using Gronwall's inequality, we derive \eqref{NN*20}. Then the remainder of the proof is as same as   in the proof of Lemma \ref{LN2}.

 (c) Proof of Lemma \ref{LJJ}. According to \cite[Lemma 2.7]{W21a}, under $(B)$ the estimate \eqref{NN9} holds for the solution to the reflecting SDE:
 $$\d \bar X_{t}=b_t^{(1)}(\bar X_t)\d t + \sigma_t(\bar X_{t}, \Phi^\gamma_{t}\mu)\d W_t+\n(\bar X_t)\d l_t\ \ \bar X_0=X_0, t\in[0,T].$$
 Then the desired assertion follows    as in the original proof.

 \end{proof}
\paragraph{Acknowledgement.} The authors   would like to thank the referee and editors for helpful comments and corrections.

\end{document}